\documentclass[12pt,twoside]{article}
\usepackage{amsfonts}
\usepackage{amsthm}
\newtheorem{theorem}{Theorem}
\newtheorem{proposition}{Proposition}
\newtheorem{conjecture}{Conjecture}
\newtheorem{lemma}{Lemma}
 
\newtheorem{corollary}{Corollary}
\newtheorem{definition}{Definition}
\newtheorem{remark}{Remark}

\def\vsni{\vskip 0.2cm}

\def\R{I\!\!R}

\def\N{I\!\!N}
\def\Q{I\!\!\!\!Q}
\def\Z{I\!\!\!\!Z}

\def\o{\omega}
\def\t{\theta}

\def\limsup{\mathop{\overline{\rm lim}}}
\def\liminf{\mathop{\underline{\rm lim}}}

\def\S{\Sigma}
\def\S2{\Sigma_2}
\def\s{\sigma}
\def\a{\alpha}

\def\ui{[0,1]}

\def\be{\begin{equation}}
\def\ee{\end{equation}}

\def\\{\hfill\break}
\def\={{ \; \equiv \; }}

\def\t_N{\tilde{\Z}_N}

\begin{document}
\date{}
%
%
%
%
\title{Recurrence and algorithmic information} 
\author{Claudio Bonanno\thanks{Dipartimento di Matematica,
Universit\`a di Pisa, via Buonarroti, 56126 Pisa, Italy. e-mail:
$<$bonanno@mail.dm.unipi.it$>$}, Stefano Galatolo\thanks{Dipartimento
di Matematica Applicata "U. Dini", Universit\`a di Pisa, via Bonanno
25/b, 56126 Pisa, Italy. e-mail: $<$galatolo@mail.dm.unipi.it$>$} and
Stefano Isola\thanks{Dipartimento di Matematica e Informatica,
Universit\`a di Camerino and INFM, via Madonna delle Carceri, 62032
Camerino, Italy. e-mail: $<$stefano.isola@unicam.it$>$.} } 
\maketitle

\begin{abstract}
In this paper we initiate a somewhat detailed investigation of the
relationships between quantitative recurrence indicators and
algorithmic complexity of orbits in weakly chaotic dynamical
systems. We mainly focus on examples.
\end{abstract}

\section{Introduction}

The general understanding of dynamical systems has recently been
enriched by two important items: the quantitative characterization of
recurrence from a statistical point of view and the measure of orbit
complexity from an algorithmic information theoretic point of view.

The basic result about recurrence is the Poincar\`e Recurrence
Theorem, which says that for a probability preserving transformation
$(X,T,\mu)$ almost all orbits starting from a set of positive measure
will come back to it infinitely many times (see, e.g., \cite{Pe}).  If
moreover $X$ is endowed with a metric $d$ then under suitable
assumptions the sequence of distances $d_n=\min_{0<i\leq n}d(T^i x,x)$
converges to $0$ as $n$ increases \cite{Fu}. This a qualitative
statement. But what is the ``typical'' speed of decreasing of $d_n$?
Boshernitzan started to investigate this problem in \cite{Bo}. His
results establish some bounds on the rate of decreasing of $d_n$ to
$0$ in terms of the Hausdorff dimension of the invariant measure. At
the same time Ornstein and Weiss proved a probabilistic formula
relating the growth of the repetition time of words of increasing
length in an ergodic source with the entropy of the source \cite{OW}.
More recently, new results have been obtained by several authors
relating various quantitative recurrence indicators to entropies and
dimensions (see, e.g., \cite{ACS}, \cite{BS},\cite{STV}).  On the one
hand, one can measure the speed of decreasing of $d_n$ or else the
increasing of the time $\tau_r$ that is necessary to come back to a
neighborhood of the starting point as the radius $r$ of the
neighbourhood decreases. Under suitable hyperbolicity assumptions this
kind of indicators are naturally related to the dimension \cite{BS}.
On the other hand, one may study the recurrence of a given orbit with
respect to the refinements of a generating partition. Under similar
hypotheses this is related to the entropy of the system (see
\cite{STV} and references therein).

Albeit the existing literaure is mainly devoted to the investigation
of systems with positive entropy, here we study recurrence properties
of some ergodic dynamical systems with zero entropy.

The entropy of a system can be interpreted as the average information
that is necessary to describe one step of the evolution of the system,
the underlying notion of information content being that due to
Shannon.  Another notion of information content is the Algorithmic
Information Content (AIC) due to Kolmogorov and Chaitin. When used to
investigate a dynamical system, the main difference comes from the
fact that the AIC is a pointwise notion.  Each digital string has its
own information content, independently form the context where the
string appears.  Therefore the corresponding amount of information
that is necessary to describe an orbit yields a pointwise version of
the entropy, called {\sl orbit complexity} in \cite{Br}.  Whenever a
system has an invariant ergodic probability measure, the orbit
complexity equals the entropy almost everywhere, hence the two notions
are quantitatively related. More precisely, if the system has positive
entropy $h>0$, the algorithmic information required to describe $n$
steps of an orbit increases as $hn+o(n)$ for almost every initial
point.

Conversely, in a null entropy system, the (sublinear) asymptotic
behaviour of the algorithmic information describing an orbit can be
related to some  dimensions and initial condition
sensitivity indicators of the system \cite{Ga}.  This suggests that
the asymptotic behaviour of the algorithmic information provides a
meaningful indicator to characterize the ``weak'' chaos present in the
system.

In the following we will see how this notion is related to
recurrence. More precisely we will see that low pointwise quantitative
recurrence rates implies strong limitations to the increase of the
information content.

The paper is organized as follows: in Section \ref{sec2} we first
recall several notions of recurrence which have been recently
investigated in the literature and also study an indicator which
measures the recurrence rate of an orbit starting from a point $x$ to
a neighbourhood of some given target point $y$ (see Theorem
\ref{Ganzo!}). This will be used in the last Section to give a precise
estimation of the information content for orbits generated by an
Interval Exchange Transformation. We then prove some general results
relating local recurrence rates to the asymptotic behaviour of the
AIC.

In a sense these results are analogous to those contained in
\cite{STV}, but they are meaningful for null entropy systems. These
results will be used in Section \ref{Markov} to estimate the local
recurrence rate (for Lebesgue almost each point) for a non-uniformly
hyperbolic system with an infinite invariant measure, for which the
standard techniques cannot be applied.

In Section \ref{sec4} we study the behaviour of the recurrence
indicators introduced in the first Section along with that of the AIC
for ergodic rotations and Interval Exchange Transformations. In
particular in Section \ref{rota} we give a precise characterization of
the recurrence rates for ergodic rotations in terms of arithmetic
properties of the angle of rotation (thus not dimension nor
entropy). From these characterizations it follows that the converse
implication ``if the information rate is low then the recurrence rate
is low'' does not hold.

In the last two Sections we study the behaviour of the AIC for
piecewise isometries using two different approaches. In Section
\ref{bona} by direct estimation of the algorithmic information
contained in the continued fraction expansion both of the initial
point and of the local rotation angle. In Section \ref{gala} by using
the recurrence rate near given target points mentioned above. This
yields sharp estimation of the information content of the orbits for a
wide class of piecewise isometries of the interval.

\section{Preliminaries and some general results} \label{sec2}

\subsection{Various kind of recurrence}\label{recu} 

Let us first recall some general facts about recurrence. Let $(X,d)$
be a complete, separable metric space and $T:X\to X$ a Borel
measurable transformation.  We shall say that a point $x\in X$ is {\rm
positively recurrent} if
\be \label{rec} 
\liminf_{n\to \infty} d(x,T^ n x) = 0.
\ee 

If $T$ is a homeomorphism the existence of such a point is ensured by
the Birkhoff Recurrence Theorem (see \cite{Pe}). More precise results
can be obtained for a metric measure preserving system (m.m.p.s)
$(X,{\cal B}, \mu, d,T)$ where the metric $d$ on $X$ is such that open
subsets of $X$ are measurable and $\mu$ is a $T$-invariant probability
measure. A theorem of Furstenberg (\cite{Fu}, p.61) asserts that for a
m.m.p.s. such that $(X,d)$ has a countable basis $\mu$-almost every
point $x\in X$ satisfies (\ref{rec}).  Quantitatively more informative
result have been obtained by Boshernitzan, Barreira and Saussol
involving various notions of dimension. For every $\alpha$, denote by
$m_\alpha (\cdot)$ the Hausdorff $\alpha$-dimensional outer measure on
$X$ (see e.g. \cite{Fa}). It is well known that for any $E \subseteq
X$ there is a unique value, ${\rm dim}_H E$, called the {\sl Hausdorff
dimension} of $E$, such that $m_\alpha (E)=\infty$ if $0\leq \alpha <
{\rm dim}_H E$ and $m_\alpha (E)=0$ if ${\rm dim}_H E<\alpha
<\infty$. The Hausdorff dimension of a probability measure $\mu$ on
$X$ is then given by
\be 
{\rm dim}_H \mu = \inf \{{\rm dim}_H E\, : \, \mu(E)=1\}.  
\ee 

Furthermore, the {\sl lower} and {\sl upper pointwise dimensions} of
$\mu$ at a point $x\in X$ are given by
\be \label{lud}
\underline{d}_\mu (x) =\liminf_{r\to 0} {\log \mu(B(x,r))\over -\log
r}\quad\hbox{and}\quad \overline{d}_\mu (x) =\limsup_{r\to 0} {\log
\mu(B(x,r))\over -\log r} 
\ee 
where $B(x,r)$ denotes the open ball of radius $r$ centered in
$x$. The measure $\mu$ is called {\sl exact dimensional} if there is a
constant $d$ such that
\be \label{exdim}
\underline{d}_\mu (x) =\overline{d}_\mu (x) = d \qquad
\hbox{$\mu$-a.e.}   
\ee 
One can show that $\underline{d}_\mu (x)\leq {\rm dim}_H \mu$ for
$\mu$-almost every $x\in X$ and in particular if $\mu$ is exact
dimensional then ${\rm dim}_H \mu =d$ (see \cite{Y}).

We now give the definitions of some used quantitative indicators of
recurrence.

First, let $s\geq 0$ be a real number and define
\be
\nu_s (x) :=\liminf_{n\to \infty} \, n^{s}\cdot d(x,T^ n x).
\ee
Given a positively recurrent point $x\in X$ let us define the sequence
$\tau_n$ of closest returns as the subsequence of the integers which
makes $d(x,T^nx)$ monotonically decreasing:
\be\label{taun} 
\tau_n(x) := \inf\, \{ k>\tau_{n-1}\, : \,
d(x,T^k x)\leq d(x,T^{\tau_{n-1}}x)\}, \quad n> 1, 
\ee
with $\tau_1$ chosen arbitrarily. Put moreover
\be\label{deltan}
d_n(x) := d\, (\, x,T^{\tau_{n}}x\, )
\ee
so that
\be\label{nus}
\nu_s (x) =\liminf_{n\to \infty} \, {\tau_n(x)}^s\cdot d_n(x).
\ee
Boshernitzan proved that given $(X,{\cal B}, \mu, d,T)$ as above if
the $\alpha$-dimensional Hausdorff measure $m_\alpha$ agrees with
$\mu$ on the $\s$-algebra ${\cal B}$ then
\be\label{bosher}
\nu_{1\over \alpha} (x) \leq 1,\qquad\hbox{$\mu$-a.e.}
\ee

If instead $m_\alpha (X)=0$ then $\nu_{1\over \alpha} (x)=0$
$\mu$-a.e. (\cite{Bo}, Thm 1.2). An improvement has been obtained in
\cite{BS}. Define the return time of a point $x\in X$ into the open
ball $B(x,r)$ by
\be\label{ret}
\tau (x,r) = \inf \{k >0 \, : \, T^k x\in B(x,r)\} = 
\inf \{k >0 \, : \, d(x,T^k x)<r\} 
\ee
(note that the sequence $\tau_n$ defined above is but
$\tau(x,d_{n-1})$). The {\sl lower} and {\sl upper pointwise
recurrence rates} of $x$ are defined as
\be\label{luprr}
{\underline {\cal R}}(x,T) = \liminf_{r\to 0 } {\log \tau(x,r)\over
-\log r} 
\quad 
\hbox{and} 
\quad
{\overline {\cal R}}(x,T) = \limsup_{r\to 0 } {\log \tau (x,r)\over
-\log r}\, \cdot 
\ee
Boshernitzan's result can then be rephrased by saying that
${\underline {\cal R}}(x,T) \leq {\rm dim}_H \mu$, $\mu$-a.e., whereas
Barreira and Saussol proved that ${\underline {\cal R}}(x,T) \leq
\underline{d}_\mu (x)$ and ${\overline {\cal R}}(x,T) \leq
\overline{d}_\mu (x)$, $\mu$-a.e. (\cite{BS}, Thm 1). In particular,
if $T$ is a ${\cal C}^{1+\eta}$ diffeomorphism acting on a hyperbolic
set and if $\mu$ is an ergodic equilibrium measure of a H\"older
continuous weight function then ${\underline {\cal R}}(x,T) =
{\overline {\cal R}}(x,T) = {\rm dim}_H \mu$, $\mu$-a.e. (\cite{BS},
Thm 5).

Let us consider a finite measurable partition ${\cal Z}$ of $X$. One
may consider cylinder sets of the form
\be\label{partset}
A=Z_n(x) \in \vee_{k=0}^{n-1}T^{-k}{\cal Z},\quad x\in Z_n(x). 
\ee
The partition ${\cal Z}$ allows us to construct a coding map $\pi : X
\to \Sigma_{\cal Z}$ onto a symbolic space $\Sigma_{\cal Z}$ which
semiconjugates $T$ with the shift $S : \Sigma_{\cal Z}\to \Sigma_{\cal
Z}$, i.e. $S \circ \pi = \pi \circ T$. Endowing $\Sigma_{\cal Z}$ with
the pseudo metric
$$
d(\o,\o') =\sum_{k\geq 0} e^{-k}\,| \o_k-\o_k'|
$$ 
we have that if $y,z\in Z_n(x)$ then $d(\pi(y),\pi(z))\leq C\,
e^{-n}$. Let us denote by $d_{\cal Z}$ the metric induced on $X$ by
$d$ via the map $\pi$, and let $h_\mu(T,{\cal Z})$ be the $\mu$-entropy
of $T$ with respect to the partition ${\cal Z}$. An easy consequence
of the Shannon-McMillan-Breiman Theorem (see \cite{CFS}, p. 257) is
that the Hausdorff dimension of $\mu$ with respect to the metric
$d_{\cal Z}$ coincides with $h_\mu(T,{\cal Z})$. Thus, in a sense, the
results mentioned above can be viewed as generalizations of a theorem
by Ornstein and Weiss (\cite{OW}, Thm 1) which in our language can be
reformulated as follows: given a measurable transformation $T:X\to X$,
a measurable partition ${\cal Z}$ of $X$ and a $T$-invariant ergodic
probability measure $\mu$ on $X$, if we endow $X$ with the metric
$d_{\cal Z}$ then ${\underline {\cal R}}(x,T) = {\overline {\cal
R}}(x,T) = h_\mu(T,{\cal Z})$, $\mu$-a.e.
\vskip 1cm

So far we considered recurrence indicators measuring how fast an orbit
comes back near to the initial point or in a set containing the
initial point.  Another kind of indicator measures how fast the orbit
of a given initial point approaches some other target point. This is
interesting when we need to understand how close the orbit of $x$
approaches some singular locus of $(X,T)$ (for example the
discontinuity set of $T$).

In the following theorem we fix a set of target points $x_i \in X$,
$i=1,\dots ,l$ and study the entrance times into the balls
$B(x_i,r)$. This result will be used in Section \ref{gala} to estimate
the type of initial condition sensitivity of a class of discontinuous
maps, the Interval Exchange Transformations.

\begin{theorem}\label{Ganzo!} 
Let $(X,{\cal B}, \mu, d,T)$ be as above and $x_1,...,x_l$ be a finite
set of points in $X$. If $\alpha>\max
\,\{{\overline{d}_\mu(x_1)}^{-1},...,{\overline{d}_\mu(x_l)}^{-1}\}$
then for $\mu$-almost all $x \in X$ it holds
$$\liminf_{n\to \infty}n^\alpha \,\min_{i}\,d(x_i,T^nx)=\infty.$$
\end{theorem}

The proof follows straightforward from the following 

\begin{lemma}\label{hjkjgh} 
Let $(X,{\cal B}, \mu, d,T)$ be as above, $x_i\in X$ and assume that
$\exists \, d$ and $\overline r >0$ such that $\mu(B(x_i,r))\leq r^d$
for each $r<\overline r$. If $\alpha>\frac 1d$ then for $\mu$-almost
every $x \in X$ we have
$$\liminf_{n\rightarrow \infty}n^\alpha \cdot d(x_i,T^nx)=\infty.$$
\end{lemma}
{\sl Proof.} Let us consider the family of sets $I_n(x_i)=
T^{-n}(B(x_i,n^{-\alpha}))$. From the fact that $\mu(B(x_i,r))\leq
r^d$ and $\alpha>\frac 1d $ it follows that $\sum _n
\mu(B(x_i,n^{-\alpha}))<\infty$.  Since $\mu$ is $T$-invariant then
$\sum _n \mu(I_n(x_i))<\infty$.  Now each point $x$ such that
$$\liminf_{n\rightarrow \infty}n^\alpha \cdot d(x_i,T^nx)< 1 $$
belongs to infinitely many sets in the family $I_n(x_i)$.
By the Borel-Cantelli Lemma this implies that $x$ 
is contained in a set of zero measure. \qed

\vsni
So long for pointwise-like results. In a slightly different perspective
we may consider a dynamical system $(X,T)$ where $T$ is a measurable
transformation along with the set $M(X)$ of $T$-invariant Borel
probability measures on $X$ and the set $EM(X)$ of ergodic
$T$-invariant Borel probability measures on $X$. Taking $\mu \in M(X)$
one can study the recurrence times into a given subset $A\subset X$ of
positive $\mu$-measure. The basic result is Poincar\'e Recurrence
Theorem which says that $\mu$-almost every point $x\in A$ returns to
$A$ infinitely often, that is $\# \, \{n>0 \, : \, T^nx\in A
\}=\infty$.

Given $A$ as above we define its {\sl Poincar\'e recurrence time}
$\tau(A)$ as
$$
\tau(A) = \min\, \{ \, n>0\,:\, T^nA\cap A \ne \emptyset \,\}=
\inf \, \{ \,\tau(x,A) \, :\, x\in A\,  \}
$$
where $\tau(x,A)= \inf\{ k>0\, : \, T^kx\in A\}$ is the return
time function in
$A$. Poincar\'e Recurrence Theorem yields for $\mu \in M(X)$
\be\label{Poi}
\mu (A)>0 \; \Longrightarrow \; \tau(A) <\infty \, .
\ee
On the other hand we have the Kac's formula (see \cite{Pe})
\be
\int_{A} \tau(x,A) \, d\mu (x) = 1.
\ee
Therefore (\ref{Poi}) can be improved to
\be
\mu (A)>0 \; \Longrightarrow \; \tau(A) \leq {1\over \mu (A)}\, \cdot
\ee
Note that Poincar\'e recurrence for a set $A$ can be very different
from typical return times $\tau(x,A)$. In the sequel we shall
consider cylinder sets of the form (\ref{partset}). In this case, as
we have seen, the return time $\tau(x,Z_n(x))$ for a $\mu$-generic
point $x$ behaves as $e^{n\, h_\mu(T,{\cal Z})}$ where $h_\mu(T,{\cal
Z})$ is the entropy of $\mu$ w.r.t the partition ${\cal Z}$. On the
other hand $\tau (Z_n(x))$ is ${\cal O}(n)$.

An elementary argument presented in \cite{STV} shows that the cylinder
$Z_n(x)$ is completely determined by its first $\tau(Z_n(x))$ symbols.
More specifically, setting $Z_n(x)=[\omega_0\omega_1\dots
\omega_{n-1}]$ and $\tau \equiv \tau(Z_n(x))$ then if $j$ denotes the
integer defined by $j\tau <n\leq (j+1)\tau$ we have
\be \label{STV}
\omega_{i\tau},\dots ,\omega_{(i+1)\tau-1}=\omega_0,\dots,
\omega_{\tau-1} \quad\hbox{for}\quad i=0,1,\dots,j\,.  
\ee

This suggests that the `effective length' of the string
$\omega_0\omega_1\dots \omega_{n-1}$ is somehow encoded by the
numerical sequence $\tau(Z_n(x))$. We may also define the {\sl lower}
and {\sl upper ${\cal Z}$-recurrence rates} as
\be\label{recrat}
{\overline R}(x,T,{\cal Z}) = \limsup_{n\to \infty }
{\tau(Z_n(x))\over n} 
\quad
\hbox{and}
\quad
{\underline R}(x,T,{\cal Z}) = \liminf_{n\to \infty }
{\tau(Z_n(x))\over n}\, \cdot\ 
\ee

We remark that this quantity is substantially different from that
defined in eq. (\ref{luprr}), that is a ratio of logarithms
to extract the power law behaviour of the recurrence times with respect
to the distance. Moreover the quantity defined by Equation
(\ref{luprr}) does not depend on the choice of a partition.
 
We now collect some results about this recurrence rate which will be
useful in what follows.

\begin{theorem}[\cite{HSV}] \label{inva} If $\mu \in M(X)$ then 
${\underline R}(x,T,{\cal Z})$ and ${\overline R}(x,T,{\cal Z})$ are
$\mu$-a.e. invariant. If moreover $\mu \in EM(X)$ then ${\underline
R}(x,T,{\cal Z})$ and ${\overline R}(x,T,{\cal Z})$ are
$\mu$-a.e. constant.
\end{theorem}

\begin{theorem}[\cite{STV}] \label{saussol} If $\mu \in EM(X)$ and
$h_\mu(T,{\cal Z})>0$ then ${\underline R}(x,T,{\cal Z})\geq 1$,
$\mu$-a.e.
\end{theorem}

\begin{remark}\label{spec}
We recall that the pair $(X,T)$ is said to have the {\rm
specification} property if the following holds: for any $\epsilon >0$,
there exists a positive integer $M(\epsilon)$ s.t. for any $k\geq 2$
and for any $k$ points $x_1,\dots , x_k \in X$, and for any positive
integers $a_1\leq b_1 <a_2\leq b_2 < \cdots a_k\leq b_k$ satisfying
$a_i-b_{i-1}\geq M(\epsilon)$ for $2\leq i \leq k$, and for any $p\geq
M(\epsilon) + b_k-a_1$, there exists a point $x\in X$ with $T^px=x$
such that $d(T^nx,T^nx_i)\leq \epsilon$ for any $a_i\leq n\leq b_i$,
$1\leq i\leq k$.  $(X,T)$ is said to be {\rm weakly specified} if the
above holds for $k=2$.  One easily sees that if $(X,T)$ is weakly
specified then $\tau(Z_n(x))\leq n$ and therefore ${\underline
R}(x,T,{\cal Z})={\overline R}(x,T,{\cal Z})= 1$, $\mu$-a.e.
\end{remark}

\begin{remark}
Theorem \ref{saussol} turns out to be a key estimate in proving
exponential and Poisson statistics for return times in dynamical
systems \cite{HSV}. In the following, however, we shall take a
different direction, by comparing the typical behaviour of
$\tau(Z_n(x))$ with the algorithmic information of the orbit of $x$
obtaining results which are meaningful in the null entropy case. These
results have the same philosophy of Theorem \ref{saussol}: if the
information rate is big then the recurrence must be big as well.
\end{remark}

\subsection{Algorithmic information content} \label{AIC}

Let us briefly recall the notion of Algorithmic Information Content
(AIC) of orbits in dynamical systems with respect to some given
partition.

To this end one starts with a Turing machine $A$ and a finite string
$s$ written in some finite alphabet $\cal A$.  Let $\Sigma$ be the set
of finite binary strings.  Let $p\in \Sigma $, $p$ will be considered
as a program to be run. If we start the machine $A$ with input $p$,
the computation stops and the output is $s$, we write $A(p)=s$. By this
notation we emphasize the function: {\em input$\rightarrow$ output}
that is naturally associated to the machine. A function is said to be
recursive if its values can be computed by a Turing machine as above.
Some computing machine stops and give an output for each given input
string.  Some other has a set of inputs that leads to a computation
that never stops.  In this case the output is not defined.

In the first case we say that the machine defines a total recursive
function (from $\Sigma$ to the set of finite strings in $\cal A$).  In
the second case to the machine it can be associated a function that is
defined only on a subset of $\Sigma$ in this case we say that $A$
defines a partial recursive function (from $\Sigma$ to another set of
strings). Let us denote by $\ell (p)$ the length of $p$.

The {\sl Kolmogorov complexity} or {\sl Algorithmic Information
Content of $s$ relative to $A$} is the quantity
$$
AIC_A(s) =\min_{p\,:\,A(p)=s}\ell (p)
$$
If there is no $p$ s.t. $A(p)=s$ then we put $AIC_A(s)=\infty$. 

There are countably many Turing machines which may be computably
enumerated as $A_1, A_2, \dots $ For $n\in \N$ let $e_n$ be the $n$-th
binary string in the lexicographical order $0, 1, 00, 01, 10, 11, 000,
\dots$ so that $\ell(e_n)\leq \log_2 n$.  We say that a Turing machine
$U$ is {\sl universal} if for any $n\in \N$ and any finite $0-1$
string $p$ we have $U({\hat e_n}\,p) = A_n (p)$ where, for a given
word $q$ of length $m$, we have denoted by ${\hat q}$ the word
$q(0)q(0)q(1)q(1)\dots q(m-1)q(m-1)01$. This in particular means that
if $A$ is any Turing machine a constant $C_A$ can be found so that for
any finite binary string $s$ we have $AIC_U(s) \leq AIC_A(s)
+C_A$. Then the AIC is independent of the choice of the universal
machine $A$ up to a constant. Since we are considering the
asymptotical behaviour of the AIC for very long strings this constant
is not relevant.

We can now define an (upper) average complexity of an infinite
sequence $\omega =\omega_0 \omega_1 \dots$ as the limit
$$
{\overline K}_U(\omega)= \limsup_{n\to \infty}{AIC_U(\omega_0\dots
\omega_{n-1})\over n}
$$
and a corresponding lower complexity ${\underline K}_U(\omega)$ with
$\limsup$ replaced by $\liminf$.

Coming back to our dynamical system $(X,T)$ and given a point $x\in X$
we may define the Algorithmic Information Content $AIC(y,n)$ of any
piece of orbit $\{y,Ty, \dots ,T^{n-1}y\}$ with $y\in Z_n(x)
=[\omega_0\omega_1\dots \omega_{n-1}]$ just as $AIC(y,n)\equiv
AIC_U(\omega_0\dots \omega_{n-1})$. In particular one defines the
upper and lower algorithmic complexities of a point $x$ w.r.t. $T$ and
${\cal Z}$ as ${\overline K}(x,T,{\cal Z})={\overline K}_U(\omega)$
and ${\underline K}(x,T,{\cal Z})={\underline K}_U(\omega)$,
respectively, where $\omega$ is the (infinite) symbolic coding of the
orbit of $x$ with $T$ on the partition ${\cal Z}$. We have the
following result

\begin{theorem}[\cite{Br},\cite{W}] \label{brudno-white} If $\mu \in
EM(X)$ then ${\underline K}(x,T,{\cal Z})={\overline K}(x,T,{\cal Z})=
h_\mu(T,{\cal Z})$, $\mu$-a.e.
\end{theorem}

The above theorem shows a strict relation between AIC and
entropy. Such a relation is also intuitive since both entropy and
$K_U$ measures the information rate of the symbolic strings that comes
from the given partition. The difference between this two notions is
that entropy is defined using the Shannon information while $K_U$ uses
the algorithmic information.  Moreover, entropy is a global (average)
notion while $K_U$ is pointwise.  The entropy and thus the increase of
the information has a crucial role in the study of chaotic dynamics
and it is related with recurrence, dimension, initial condition
sensitivity.  By the above theorem also the algorithmic information of
symbolic orbits is related to all these indicators of chaos.

In the case where entropy is null the information that is necessary to
describe the orbit increases less than linearly. There are many
possible such asymptotic behaviours. In this case then the asymptotic
behaviour of the algorithmic information that is necessary to describe
the orbit gives an orbit complexity indicator more refined than $K_U$
(that is almost always null). This information behaviour is also
related to dimension and initial condition sensitivity. Recent
results in \cite{Ga} proves quantitative relations between dimension,
AIC and initial conditions sensitivity that are meaningful in the zero
entropy case.

In the following we start to investigate also some relation with
recurrence. We will see that certain recurrence rates implies strong
limitations to the algorithmic information of the orbits. In Section
\ref{gala} we will also see how by the recurrence near given points it
is possible to have an estimation of the initial condition sensitivity
of a system, and then of the algorithmic information of its orbits.

\vsni
\vsni

As noted above (Eq.(\ref{STV})) the AIC of any $n$-long orbit of
points in $Z_n(x) =[\omega_0\omega_1\dots \omega_{n-1}]$ is bounded by
the AIC of the initial word $\omega_0\dots \omega_{\tau-1}$ plus the
AIC needed to repeat this word up to the size $n$ of the cylinder,
that is 
\be\label{basic} 
AIC(x,n)\leq AIC(x,\tau(Z_n(x))) + \log n +C
\ee

To construct the symbolic string that gives the $n$ steps symbolic
orbit of $x$ with respect to ${\cal Z}$ we need the information that
is necessary for the first period (that is $\tau(Z_n(x))$ steps long)
and the information needed for $n$ (that is about $\log n $ bits).

If, moreover, $(X,T)$ is weakly specified then $\tau(Z_n(x))\leq n$ so
that we also have
\be
AIC(x,\tau(Z_n(x)))\geq AIC(x,n).
\ee

Therefore, in this case the condition $h_\mu(T,{\cal Z})>0$
entails\footnote{{\sc Notations:} Here and in the sequel, for two
sequences $a_n$ and $b_n$ we shall write $a_n \sim b_n$ if the
quotient $a_n/b_n$ tends to unity as $n\to \infty$.  Moreover, the
notation $a_n\approx b_n$ means that $a_n/b_n=o(n^{\epsilon})$ as well
as $b_n/a_n=o(n^{\epsilon})$ for $n\to \infty$ and $\forall \epsilon
>0$. This condition is satisfied if $a_n/b_n \sim L(n)$ where $L(n)$
is some function slowly varying at infinity, i.e. $L(cn) \sim L(n)$
for every positive $c$.  Both $\sim$ and $\approx$ are equivalence
relations and we shall denote by $[a_n]$ and $[[a_n]]$ respectively,
the equivalence classes of $a_n$ w.r.t to these relations.}  
\be
AIC(x,n) \sim h_\mu(T,{\cal Z})\ \tau(Z_n(x)) \sim h_\mu(T,{\cal Z})\
n,\qquad \mu-a.e.  
\ee

On the other hand, by Theorem \ref{brudno-white} the condition
$h_\mu(T,{\cal Z})=0$ implies ${\underline K}(x,T,{\cal Z})={\overline
K}(x,T,{\cal Z})=0$, $\mu$-a.e. and therefore 
\be 
AIC(x,n) = o(n),\qquad \mu-a.e.  
\ee 

Can we say something about the asymptotic relations between $AIC(x,n)$
and $\tau(Z_n(x))$ in this case?  We start proving the following
results that are in some sense extensions of the Theorem \ref{saussol}
for the null entropy case. We include the case of systems with an
infinite invariant measure, assuming conservativity of the systems, a
classical assumption in these cases that is verified for systems with
a probability invariant measure. The extension to infinite measures
will prove useful in Section \ref{Markov}.

\begin{theorem} \label{aic} Let $(X,T,\mu)$ be an ergodic conservative
dynamical system, with $\mu$ a $T$-invariant and not necessarily
finite measure on $X$. Then
$$
\bar \alpha (x) := \inf \left\{ \alpha \ : \  \liminf_{n\to \infty} \
\frac{\tau(Z_n(x))}{n^\alpha}  = 0 \right\} 
$$
is $\mu$-a.e. equal to a constant $\bar \alpha$. If moreover
$0\leq {\overline K}(x,T,{\cal Z}) <\infty$ 
then for $\mu$-a.e. $x\in X$ we have
$$\liminf_{n\to \infty} \frac{AIC(x,n)}{n^\beta} =0,\qquad \forall \, \beta > \bar \alpha .$$
\end{theorem}
{\sl Proof.} The constancy of $\bar \alpha (x)$ $\mu$-a.e. follows at
once from Theorem \ref{inva}. Using (\ref{basic}) we have, for $\beta
> \bar\alpha$,
\begin{eqnarray}
0&\leq & \liminf_{n\to \infty} \frac{AIC(x,n)}{n^{\beta}} \le
\liminf_{n\to \infty} \left( \frac{AIC(x,\tau(Z_n(x)))}{n^{\beta}} +
\frac{\log n}{n^{\beta}} \right) \nonumber \\ &\le& \limsup_{n\to
\infty} \frac{AIC(x,\tau(Z_n(x)))}{\tau(Z_n(x))} \cdot \liminf_{n\to
\infty} \frac{\tau(Z_n(x))}{n^{\beta}} = {\overline K}(x,T,{\cal Z})
\cdot 0 =0 \nonumber
\end{eqnarray}
where we recall that if $\mu(X)<\infty$ then by Theorem
\ref{brudno-white} it holds $h_\mu (T,{\cal Z}) = {\overline
K}(x,T,{\cal Z})$, $\mu$-a.e. If instead $\mu(X)=+\infty$ then
$\mu$-a.e. it holds ${\overline K}(x,T,{\cal Z})=0$ (see \cite{B}
Thm. 3.3). \qed

\vsni
From this theorem we may deduce some consequences. 

\begin{corollary} \label{h=0} Let $\mu \in EM(X)$ and suppose that for
$\mu$-a.e. $x\in X$ we have $\bar \alpha (x)=\bar \alpha< 1$. Then
$h_\mu(T,{\cal Z})= 0$. If moreover $\bar \alpha =0$ then
$AIC(x,n)\approx {\rm constant}$ $\mu$-a.e.
\end{corollary}

We remark that the first part of this statement is contained in
Theorem \ref{saussol}.

If $\overline{ R} (x,T,{\cal Z})<1$ the system is forced to have a
strong constraint about the information content of its orbits, as it
is shown in the following theorem:

\begin{theorem}\label{overR}
If $\overline{ R} (x,T,{\cal Z})<1$ then $AIC(x,n)\approx {\rm
constant}$. More specifically we have
$$\limsup_{n\to
\infty}
\frac{{AIC}(x,n)}{\log^2n}\leq{\frac {-1}{\log \overline {R}(x)}}.$$
\end{theorem}

The proof of this theorem consists in a straightforward application
to Eq.\ref{basic}   of
the following lemma with the identifications $d=1 /{\overline R}$ and
$a_n=AIC(x,n)$.

\begin{lemma}\label{mos}
Let $a_n$ be a real sequence such that $\exists d>1$ and $\overline n$
such that $\forall n >\overline n$ there is $d_n $ such that $d_n\geq
d$, $ \frac n{d_n}\in \N$ and $a_n\leq a_{\frac n {d_n}} +\log n+
C$. Then
$$\limsup_{n\to \infty} \frac {a_n}{\log^2 n}\leq\frac 1{\log d}$$
\end{lemma}
{\sl Proof.} Set $A=\max_{n\leq {\overline n}}\, a_n$ and $f:\N \to
\N$ be the map given by $f(n)={n\over d_n}$. Note that $f(n) <n$.
Letting $n>\overline n$ we have
$$
a_n\leq a_{f(n)} +\log n +C \leq a_{f^2(n)} +\log n+\log f(n)+
2C\leq \cdots \leq a_{f^k(n)}+k\, (C+\log n)$$ with $k=\min \{i\in \N : f^i(n)\leq { \overline n} \}$. One
readily realizes that
$k\leq
\frac {\log\, ( n/{ \overline n} )}{\log d}$.
Therefore
$a_n\leq A+{\log\, ( n/ {\overline n})\over \log d}\, (C+\log n)$. \qed 

\vsni 
We remark that the above result is pointwise and does not depend on
the invariant measure. We will see an application of Theorems
\ref{aic} and \ref{overR} in the next section for dynamical systems
with an infinite invariant measure absolutely continuous with respect
to the Lebesgue measure.

Of course Theorem \ref{aic} is interesting when $\bar \alpha <
1$. When $\bar \alpha =1$ the conclusion of Theorem \ref{aic} is not
so interesting. However, as we shall see below (Theorem \ref{bounded})
there are uniquely ergodic dynamical systems with zero entropy and
such that $0<{\underline R}(x,T,{\cal Z})\leq {\overline R}(x,T,{\cal
Z}) \leq 1$ (and thus $\bar \alpha \geq 1$). By the way, this provides
a counterexample to show that the inequality in the above theorem
cannot be reversed in general. There are systems with very low orbit
complexity and several possible recurrence rates. We are now going to
discuss some examples.

\section{Markov interval maps}\label{Markov}

In this Section we apply the results of the previous one to
calculate recurrence times for a class of intermittent maps on the
interval, with finite or infinite absolutely continuous invariant
measure.

Roughly speaking, a map $T:[0,1]\to [0,1]$ is called Markov if
there is a finite partition ${\cal Z}$ of $[0,1]$ such that the image
of each of its elements is the union of some of its elements.

Let us consider a map $T$ such that 

1) $T$ is increasing and ${\cal C}^{2}$ restricted to each interval of
   the partition ${\cal Z}$;

2) $T$ is onto on each interval of the partition.

All maps satisfying 1) and 2) are obviously Markov. We remark that all
the maps that have properties 1) and 2) are topologically conjugate,
but the conjugation may be not absolutely continuous and then some
properties which are related to the Lebesgue measure may change.

Let us consider the partition ${\cal Z}=\{[0,1/2),[1/2,1]\}$. The map
$\pi : [0,1]\to \{0,1\}^{\N}$ defined by $\pi(x)
=\omega_0\omega_1\omega_2\dots$ according to $T^jx\in I_{\o_j}$ for
$j\geq 0$ is a coding map which is a homeomorphism on the residual set
of points in $[0,1]$ which are not preimages of $1$.  Moreover $S\circ
\pi = \pi \circ T$ where $S$ denotes the right shift map on
$\{0,1\}^{\N}$.  It is plain that there are exactly $2^n$ cylinders of
the form $Z_n(x)=[\omega_0,\omega_1,\dots, \omega_{n-1}]$ and each of
them contains exactly one periodic point of prime period $n$. This
implies at once $\tau(Z_n(x))\leq n$ for every $x$.  Now, if $T$ is
uniformly expanding, i.e. $\exists \rho >1$ such that $T'(x) \geq
\rho$ for all $x\in [0,1]$ then it is well known (see \cite{KH}) that
$T$ possesses a unique (ergodic) absolutely continuous invariant
measure $\mu$ with $h_\mu(T)=h_\mu(T,{\cal Z})>0$. Putting together
the above, Theorem \ref{saussol} and Remark \ref{spec}, we get
${\underline R}(x,T,{\cal Z})={\overline R}(x,T,{\cal Z})= 1$,
$\mu-a.e$. This condition continues to hold for the intermittent
situation where $T'(0)=1$ provided the invariant measure $\mu$ is
finite (so that the entropy is still positive).

In the opposite case (i.e. for example when $T'(0)=1$ and $T''(0)=0$) we have no absolutely continuous invariant probability measure.
In such case in principle we only have the inequalities ${\underline
R}(x,T,{\cal Z})\leq {\overline R}(x,T,{\cal Z})\leq 1$, $\mu-a.e.$.
In this case however, if the map satisfies some weak auxiliary
assumptions, ${\overline R}$ can be estimated.

\begin{proposition} \label{mann1}
Under the assumptions 1) and 2), if $T$ is such that $T'(0)=1$ and
$T''(x)\sim x^{z-2}$ as $x\to 0$, for some $z> 2$, and there exists an
absolutely continuous infinite invariant measure, then ${\overline
R}(x,T,{\cal Z})= 1$ $Lebesgue-a.e.$
\end{proposition}
{\em Proof.} First, we can use the results of \cite{Pr} to find out a
$C^1$-diffeomorphism of $\ui$ which conjugates $T$ with the map
$T_z(x)=x+x^z \ ({\rm mod} 1) $. In turn, the latter (called
Manneville map) is absolutely continuously conjugated to its piecewise
linear versions that can be studied by the theory of recurrent events
(see \cite{BoGa}), obtaining that $AIC(x,n)$, with respect to the
partition $\cal Z$, increases for $z>2$ as a power law ${n^\alpha}$,
with $\alpha=\frac{1}{z-1} \in (0,1)$, and that for Lebesgue a.e. $x$
it holds
\begin{equation} \label{liminfmann}
\liminf_{n\to \infty} \ \frac{AIC(x,n)}{n^\beta} = +\infty \qquad
\forall \ \beta < \alpha
\end{equation}

Hence if ${\overline R}(x,T,{\cal Z})< 1$ on a set of positive
Lebesgue measure, then by Theorem \ref{overR}, $AIC(x,n)\approx {\rm
constant}$ on this set, and this would contradict Equation
(\ref{liminfmann}). \qed

\vsni
An analogous result can be obtained for the index $\bar \alpha$ of
Theorem \ref{aic}.

\begin{proposition} \label{mann2} Under the same hypothesis of
Proposition \ref{mann1} it holds $\bar \alpha \ge \frac{1}{z-1}$.
\end{proposition}
{\sl Proof.} By contradiction, if $\alpha < \frac{1}{z-1}$, then by
Theorem \ref{aic} it follows that Equation (\ref{liminfmann}) cannot
be verified. \qed

\section{Ergodic piecewise isometries} \label{sec4}

\subsection{Recurrence results for ergodic rotations}\label{rota}

Let $X=(0,1]$ be the unit circle and $\mu$ the Lebesgue measure. Let
moreover $\alpha \in (0,1]$ and $T: X\to X$ be the rigid translation
\be
T (x) = x+\alpha \;\; ({\rm mod}\, 1).
\ee
The measure $\mu$ is clearly $T$-invariant. If $\alpha \in \R\setminus
\Q$ all orbits are dense and $\mu$ is the unique (and thus ergodic)
invariant measure, which is obviously exact dimensional with $d=1$
(cfr. (\ref{exdim})).  Let us consider the partition ${\cal Z} =
(A_0,A_1)$ of $X$ into the half-open arcs $A_0=[0, 1/2)$ and
$A_1=[1/2,1)$.  For $\alpha$ irrational $T$ is a minimal homeomorphism
of $X$. Therefore ${\cal Z}$ is generating for $T$ and the partition
${\cal Z}^{n} ={\cal Z} \vee T^{-1}{\cal Z} \vee \cdots \vee
T^{-n+1}{\cal Z}$ is made out of $2n$ arcs. This can be easily
realized by induction: when passing from ${\cal Z}^{n-1}$ to ${\cal
Z}^{n}$ one has to add to the endpoints of the arcs belonging to
${\cal Z}^{n-1}$ the two new points $T^{-n}(0)$ and
$T^{-n}(1/2)$. Thus,
$$
0\leq  h_\mu(T,{\cal Z})\leq \lim_{n\to \infty} {\log (2n)\over n} =0.
$$
It is interesting to note that the lengths of the arcs of ${\cal
Z}^{n}$ can actually take only five values (see \cite{AB}, Thm 18).

Note moreover that $(X,T)$ is metrically isomorphic to the subshift
given by the closure of $\pi ([0,1))$ where the coding map $\pi :
[0,1]\to \{0,1\}^{\N}$ is given by $\pi(x)_n=\chi_{A_1} (T^n(x))$.
Let
\be\label{cfe}
\alpha = {1\over \displaystyle a_1 + {1\over \displaystyle a_2 +
{1\over\displaystyle a_3 +\cdots }}}\equiv [a_1,a_2,a_3,\dots]
\ee 
be the continued fraction expansion of $\alpha$ and $p_n/q_n
=[a_1,a_2,\dots, a_n]$ the sequence of its rational approximants. It
is well known that (see \cite{RS})
\be\label{estima}
{1\over q_n\cdot \, (q_{n+1}+q_n)} < \left|\alpha - {p_n\over q_n}\right| < 
{1\over q_n\cdot \, q_{n +1}}.
\ee
From the estimate (\ref{estima}) one deduces the following
characterization for the denominators $q_n$: let $\Vert x\Vert=\min
\{|x-p\,|\, : \, p\in \Z\}$ be the distance from the nearest integer,
then
\be\label{magg}
q_{n} = \min \, \{ k > q_{n-1} \, : \, \Vert k  \alpha  \Vert < \Vert
q_{n-1} \alpha  \Vert \, 
\}
\ee
so that $\Vert k \alpha \Vert \geq \Vert q_{n-1} \alpha \Vert$ for all
$k < q_{n}$. Also notice that
\be\label{dist0}
d(x,T^k x)=\min_{p\in \Z}|x+k \alpha -x-p|= \Vert k
\alpha\Vert 
\ee
for all $x\in X$. Now set 
\be
f_n:= d(x,T^{q_n}x)=\Vert q_n \alpha\Vert=|q_n \alpha -p_n|=(-1)^n(q_n
\alpha -p_n) 
\ee
It holds
\be\label{recursio}
f_{n-1}= a_{n+1}f_{n}+f_{n+1}
\ee
and also
\be\label{decr}
f_n =\prod_{k=0}^{n}G^{k}(\alpha)
\ee
where $G :[0,1]\to [0,1]$ is the Gauss transformation given by $G(x)
=\{1/ x\}$ for $x> 0$ and $G(0)=0$. Equation (\ref{decr}) shows that
for all irrational $0<\alpha <1$ the sequence $f_n$ is strictly
decreasing and by (\ref{estima}) it satisfies the bounds
\be\label{dist}
{1\over a_{n+1}+2}<q_n\cdot f_n  <{1\over  a_{n+1}}
\ee

\begin{definition}
We shall say that the irrational number $\alpha$ is of type $\gamma$
if 
\be \label{nu} 
\gamma = \sup \ \left\{ \beta >0 \ / \ \liminf_{k\to \infty} k^\beta
\cdot \Vert k\, \a \Vert =0 \right\} .    
\ee
\end{definition}

By (\ref{magg}) one has that if $q_n\leq k <q_{n+1}$ then
$k^\beta\cdot \Vert k\, \a \Vert \geq q_n^\beta \cdot \Vert q_n\, \a
\Vert$ and the lower limit in (\ref{nu}) is reached along the
subsequence $k=q_n$. In addition by (\ref{dist}) we have
$q_n^{1-\epsilon} \cdot \Vert q_n\, \a \Vert < 1/q_n^\epsilon$ and
therefore $\gamma \geq 1$.  The set of type $1$ numbers includes those
with bounded partial quotients (see below, Proposition
\ref{bounded}). One has (see \cite{RS}, cap. 4.5-4.6) $\liminf_{k\to
\infty} k \cdot \Vert k\, \a \Vert \leq 1/\sqrt{m^2+1}$ for any
$\a=[a_1,a_{2},\dots]$ unless $a_i <m$ for all $i$ large enough. For
the golden mean $\alpha =[1,1,1,\dots] = {\sqrt{5}-1\over 2}$ we find
$\liminf_{k\to \infty} k \cdot \Vert k\, \a \Vert=1/\sqrt{5}$.  On the
contrary, we say that $\alpha$ is well approximated by rational
numbers if $\eta (\alpha)>0$ where $\eta (\alpha)$ is the supremum of
all $\eta$ such that for all $k\geq 1$ we have $\Vert q_k\, \a \Vert
\leq C\, q_k^{-(1+\eta)}$. In this case $q_k \cdot \Vert q_k \, \a
\Vert \leq C\, q_k^{-\eta (\alpha)}$ and $\alpha$ is of type $1+\eta
(\alpha)$. If for example $a_k=2^{2^k}$ then $\alpha$ is of type $2$.

Let us first dwell upon pointwise recurrence properties. First,
Boshernitzan theorem mentioned in Section \ref{recu} along with
minimality yield that the function $\nu_1(x)$ defined in (\ref{nus})
is finite for all $x$ and all $\alpha \in \R\setminus \Q$, which is
not very informative. We are now going to show that, in some sense,
for ergodic rotations the role of the Hausdorff dimension is played by
the type of the irrational number $\alpha$. First, from the general
properties of continued fractions recalled above one readily gets the
following

\begin{theorem} Let $X=\R/\Z$ and $T:X\to X$ be the rotation of an
irrational angle $0<\alpha <1$. Then the quantities $\tau_n$,
$d_n$ and $\nu_s$, defined in (\ref{taun}), (\ref{deltan}) and
(\ref{nus}), respectively, are constant on $X$ with $\tau_n = q_n$ and
$d_n = f_n$. Moreover, if the irrational number $\alpha$ is of
type $\gamma \geq 1$ then $\nu_\gamma <\infty$ and $\nu_s =0$ for
$s<\gamma$.
\end{theorem}

In terms of the lower pointwise recurrence rate ${\underline {\cal
R}}(x,T)$ defined in (\ref{luprr}), the above result is equivalent to
the following result

\begin{theorem} If $\alpha$ is of type $\gamma \geq 1$ then
${\underline {\cal R}}(x,T)= 1/\gamma$. 
\end{theorem}
{\sl Proof.}  Note that in this case the return time defined in
(\ref{ret}) can be written as $\tau (x,r) =\inf \{ k>0 : \Vert k\, \a
\Vert < r\}$. Recalling that the numbers $f_n$ are strictly
decreasing, for each $r$ let $n$ be such that $f_n< r \leq
f_{n-1}$. By (\ref{magg}) we have that for $0<k<q_n$ it holds $\Vert
k\, \a \Vert \geq \Vert q_{n-1}\, \a \Vert =f_{n-1} \geq r$ but $\Vert
q_{n}\, \a \Vert =f_{n}<r$. Therefore, for these values of $r$ we have
$\tau (x,r)= q_n$.  If moreover $\alpha$ is of type $\gamma$ then for
$\epsilon >0$ and $k>0$ there is $C_\epsilon>0$ s.t. $k^{\gamma
+\epsilon}\cdot \Vert k\, \a \Vert \geq C_\epsilon$. In particular
choosing $k=q_n$ this gives
$$
q_n^{\gamma +\epsilon}> {C_\epsilon\over f_n} \geq {C_\epsilon\over
r}\, \cdot 
$$
Therefore $(\gamma +\epsilon)\log q_n \geq -\log r + \log C_\epsilon$
and
$$
{\underline {\cal R}}(x,T)=\liminf_{r\to 0} {\log (\tau (x,r)) \over
-\log r}=\liminf_{r\to 0} {\log q_n \over -\log r}\geq {1\over \gamma
+\epsilon},  
$$
whence ${\underline {\cal R}}(x,T) \geq 1/\gamma$ by the arbitrariness
of $\epsilon$.

In the other direction note that $\liminf k^{\gamma-\epsilon}\cdot
\Vert k\, \a \Vert =0$. Therefore there is a subsequence $k_j$ tending
to infinity so that $k_j^{\gamma-\epsilon}\cdot \Vert k_j\, \a \Vert
<1$. Choose a sequence $r_j$ of $r$-values such that $r_j^{-1}<
k_j^{\gamma-\epsilon} \leq 2\, r_j^{-1}$. It then follows that $\Vert
k_j\, \a \Vert < r_j$ and therefore $\tau (x, r_j) \leq k_j \leq
(2/r_j)^{1\over \gamma-\epsilon}$. Finally we get
$$
{\log \tau(x,r_j)\over -\log r_j }\leq \left({\log 2 -\log r_j \over
-\log r_j}\right)\, {1\over \gamma-\epsilon}\,  
$$
Taking the limit $j\to \infty$ we get ${\underline {\cal R}}(x,T) \leq
1/\gamma$ by the arbitrariness of $\epsilon$. \qed

\vsni 
As far as the recurrence rates (\ref{recrat}) are concerned we first
mention a result obtained by Afraimovich, Chazottes and Saussol which
in our language says that if $\alpha$ is of type $\gamma >3$ then
${\underline R}(x,T,{\cal Z})=0$ (see \cite{ACS}, Thm 4.4).

\vsni
Now, according to (\ref{dist0}) the distance $d(x,T^n x)$ does not
depend on $x$ but only on $n$ and $\alpha$, so that for a given
(half-open) arc $A\subset X$ of length $|A|$ we have
\be\label{tau}
\tau (A) = \min \{ r>0 \, : \;  \Vert r\alpha\Vert < |A|\}.
\ee
Therefore for two half-open arcs $A$ and $B$ we have
\be\label{<}
|A| < |B|  \Longrightarrow \tau (A)  \geq \tau (B)
\ee
A more precise information can be obtained from the {\sl three gap
theorem} (see \cite{AB}, Sec. 4). First, it says that if $A$ is such
that $|A| \leq 1/2$ then $|A|$ can be expressed uniquely as
\be\label{rapp}
|A| = c\, f_k + f_{k+1} + g\quad\hbox{{\rm for}{  } {\rm some}}\quad
k\geq 1, \quad 0< g \leq f_k, \quad 1\leq c \leq a_{k+1}.
\ee
Using (\ref{recursio}) we see that the integer $k$ is determined by
the inequalities
$$
f_k + f_{k+1} <|A| \leq f_{k-1} + f_{k} 
$$
Second, the return time to the set $A$ can assume, on the points of
$A$, at most three values, one being the sum of the other two: they
are $r_1=q_k$ (with frequency $(c-1)f_k+f_{k+1}+g$), $r_2=q_{k+1}-c\,
q_k$ (with frequency $g$) and $r_3=r_1+r_2$ (with frequency $f_k-g$).

This yields the following
\begin{lemma}\label{dico} Let $|A|$ be as in (\ref{rapp}). Then
$$
\tau (A)=\cases{q_k , &if $\; c<a_{k+1}$, \cr q_{k-1}, &if $\;
c=a_{k+1}$. \cr }
$$
\end{lemma}
{\sc Example.} If $|A|=f_{k-1}$ for some $k>1$ then using
(\ref{recursio}) we find $g=f_k$, $c=a_{k+1}-1$ and $\tau
(A)=q_{k}$. Note that the same results from (\ref{tau}) and
(\ref{magg}).  In this case the largest return time $r_3$ has
frequency zero and therefore, by uniform distribution, it does not
appear at all.

\begin{theorem}\label{bounded}
Suppose that $\alpha =[a_1,a_2,\dots ]$ with $a_i=O(1)$, $\forall i$,
then ${\underline R}(x,T,{\cal Z})>0$ for all $x$.
\end{theorem}
{\sl Proof.} By the {\sl five distance theorem} (see \cite{AB},
Sec. 7) we know that the lengths of the intervals forming the
partition ${\cal Z}^{n}$ take at most five values. In order to
estimate the lengths of such intervals let us reason as follows.  Let
${\cal P}^n$ and ${\cal Q}^n$ be the partitions of the circle by the
points $\{j\alpha\}_{0\leq j <n}$ and $\{{1\over 2}+j\alpha\}_{0\leq j
<n}$, respectively. Then we have ${\cal Z}^{n}={\cal P}^n\vee {\cal
Q}^n$. From the {\sl three distance theorem} (see \cite{BCF}, Sec. 3)
we know that for $q_{k}<n\leq q_{k+1}$ the intervals of the partitions
${\cal P}^{n}$ and ${\cal Q}^{n}$ have three possible lengths, which
are given by: $f_{k}$, $c_n\, f_{k}+f_{k+1}$, for some $1\leq c_n \leq
a_{k+2}$, and $(1+c_n)\, f_{k}+f_{k+1}$ (this last one being not taken
for some values of $n$, in particular for $n=q_{k+1}$).  Therefore,
for $q_{k}<n\leq q_{k+1}$ the intervals of ${\cal Z}^{n}$ have lengths
not larger than $(1+c_n)\,f_{k}+f_{k+1}$.

Now, if $a_i \leq m$ then for $q_{k}<n\leq q_{k+1}$ the intervals of
${\cal Z}^{n}$ have lengths not larger than
$(m+1)\,f_{k}+f_{k+1}$. Let $B$ be such that $|B|= m\, f_{k} + f_{k+1}
+ g$ with $0<g\leq f_{k}$. By the above, we can find $k'\leq k$ with
$k-k'$ bounded uniformly in $k$, together with $0< g' \leq f_{k'}$ and
$1\leq c \leq a_{k'+1}$ so that $|B| = c\, f_{k'} + f_{k'+1} +
g'$. Whence $\tau (B)\geq q_{k'-1}$ by Lemma \ref{dico}. Moreover,
using (\ref{<}) we see that $\tau (A)\geq q_{k'-1}$ for any $A\in
{\cal Z}^{n}$ with $q_{k}<n\leq q_{k+1}$. This yields
$$
{\underline R} (x,T,{\cal Z})\geq \liminf_{k\to \infty }
{q_{k'-1}\over q_{k+1}}>\liminf_{k\to \infty }
\prod_{j=k'}^{k+1}(a_{j}+1)^{-1} \geq (m+1)^{k'-k-1}> 0
$$
where we have repeatedly used the inequality
$q_{k+1}<(a_{k+1}+1)q_k$. \qed

\begin{remark} The above result can be further improved for some
special cases. Let us consider for example $\alpha
=(\sqrt{5}-1)/2\simeq 0.62$.  $\alpha$ is the positive root of the
equation ${1\over \alpha}=1+\alpha$ and therefore its partial
quotients are all equal to one: $a_n=1$, $\forall n$. Hence for
$q_{k}<n\leq q_{k+1}$ the intervals of ${\cal Z}^{n}$ have lengths not
larger than $2\,f_{k}+f_{k+1}$ ($m=1$ and $k'=k$ in the proof above).
Thus
$$
{\underline R} (x,T,{\cal Z})\geq \lim_{k\to \infty } {q_{k-1}\over
q_{k+1}}=\alpha^2={3-\sqrt{5}\over 2},
$$
where for the last identity we have used the fact that in this example
the $q_k$'s are the Fibonacci numbers $F_k$ defined as $F_{0}=F_1=1$,
$F_k=F_{k-1}+F_{k-2}$, $k\geq 2$, and satisfying $F_k = (\alpha^{-k}
-(-\alpha)^{k})/ \sqrt{5}$.
\end{remark}
We can now formulate the following\footnote{When this paper was almost
finished we received a preprint by M Kupsa entitled {\sl Local return
rates in Sturmian shifts}, whose techniques and results seem promising
to prove our Conjecture.}

\begin{conjecture} For all irrational $\alpha$ we have $0\leq
{\underline R}(x,T,{\cal Z}) <1$, $\mu$-a.e. Moreover ${\underline
R}(x,T,{\cal Z}) >0$ $\mu$-a.e. if and only if $\alpha$ is of type
$1$.
\end{conjecture}

\subsection{Interval Exchange Transformations} \label{iet}

A rotation of the circle is the simplest example of an Interval
Exchange Transformation (IET) defined as follows (see \cite{CFS} for
details). Let $X=[0,1)$ as above and let $\xi=\{I_1,\dots ,I_r\}$,
$r\geq 2$, be a partition of $X$ into half-open intervals, numbered
from left to right, and $\sigma =(\sigma_1,\dots, \sigma_r)$ a
permutation of $(1,2,\dots ,r)$. On each half-interval $I_i$ the map
$T$ acts as a translation $T_{\alpha_i}x=x+\alpha_i$ in such a way
that the intervals get `exchanged' according to the permutation
$\sigma$. The new intervals $TI_i=T_{\alpha_i}I_i=I_i'$ are thus
ordered as $I'_{\sigma_1},\dots ,I'_{\sigma_r}$. In particular the
angles $\alpha_1, \dots, \alpha_r$ are uniquely determined by the pair
$(\xi, \sigma)$. More precisely it holds
\be\label{angles}
\alpha_i=\sum_{k<\sigma^{-1}(i)}|I_{\sigma_k}| - \sum_{k<i}|I_{k}|. 
\ee

\vsni
\noindent
{\sc Example.} Let $a,b$ be two real numbers such that $0<a<b<1$ and
consider the partition of $[0,1)$ into the half-intervals $I_1=[0,a)$,
$I_2=[a,b)$ and $I_3=[b,1)$, along with the permutation $\sigma
=(3,2,1)$. Then we have $I'_{\sigma_1}=I'_3 =[0,1-b)$,
$I'_{\sigma_2}=I'_2 =[1-b,1-a)$ and $I'_{\sigma_3}=I'_1
=[1-a,1)$. Therefore $T$ is defined as 
\be 
T x= \cases{ x+1-a &if $x\in I_1$, \cr 
x+1-(a+b) &if $x\in I_2$,\cr 
x-b& if $x\in I_3$.\cr} 
\ee 
Note that $T$ is continuous everywhere but in $a$ and $b$.  

\vsni 
It is easy to see that either $T$ has no periodic points or there
exists an interval entirely formed by periodic points of the same
period. In the former case $T$ is said {\sl aperiodic}. If we denote
by $D$ the set of discontinuity points for $T$, that is the set of
left-endpoints $d_i$, $i=1,\dots ,r$, of the intervals exchanged by
$T$, then aperiodicity is equivalent to the fact that
$\cup_{k=0}^{\infty}T^{k}D$ is dense in $[0,1)$.  Furthermore, if the
full orbits of the points $d_i$ are mutually disjoint infinite sets
then the dynamical system $(X,T)$ is minimal and uniquely ergodic
\cite{CFS}.  Finally, an interval exchange transformation cannot be
strongly mixing, but it can be weakly mixing \cite{FHZ}.

As far as recurrence properties are concerned, very few quantitative
result are available. Boshernitzan proved that for an aperiodic
exchange of $r$ intervals such that the lengths of all exchanged
intervals are algebraic numbers the Hausdorff dimension of any
$T$-invariant continuous probability measure $\mu$ on $X$ is bounded
below by $1/(r-1)$ (see \cite{Bo}, Corollary 6.9). Therefore, by
(\ref{bosher}) we have that $\nu_{s} (x) \leq 1$, $\mu$-a.e., for some
$s\leq r-1$. On the other hand this result can be interesting only for
$r\geq 4$, since for $r=2,3$ it is known that $\mu$ must be the
Lebesgue measure and thus $s=1$. In view of the results discussed in
the previous subsection we argue that in this case as well recurrence
behaviour should be somehow related to arithmetic properties of the
lengths of the exchanged intervals or alternatively, by
(\ref{angles}), to those of the local rotation angles $\alpha_i$.

\subsection{AIC for rotations and Interval Exchange
Transformations} \label{bona}

We first study the asymptotic behaviour of the AIC for the symbolic
sequences generated by the ergodic rotation $(X,T)$, with $X=(0,1]$,
$T(x)=x+\alpha \ (\mbox{mod }1)$ and $\alpha \in \R
\setminus \Q$. To study the AIC for IET we will repeatedly apply the same
argument.

When studying the behaviour of the AIC, it turns out to be important
the knowledge of the continued fraction expansion of both the number
$\alpha$ and of the initial condition. We shall use the following
estimate for the approximation of a number $\alpha$ by the rationals
$\frac{p_k}{q_k}$, following at once from (\ref{estima}):
\be \label{eq:approx_cf} 
\frac{1}{q_k^2\cdot \, (a_{k+1}+2)} <
\left|\alpha - \frac{p_k}{q_k}\right| < \frac{1}{q_k^2\cdot \,
a_{k+1}}
\end{equation}

Let $\sigma$ denote the infinite word in $\N^{\N}$ associated to the
continued fractions expansion $[a_1,a_2,\dots]$ of the angle of
rotation $\alpha$. In the same way let
$[\xi_1,\xi_2,\dots,\xi_n,\dots]$ be the continued fraction expansion
for a number $x\in (0,1)$ (the initial condition), and denote by $s$
the associated word. We also denote by $\frac{p_k}{q_k}$ and
$\frac{r_h}{t_h}$, respectively, the rational approximants of $\alpha$
and $x$:
$$\frac{p_k}{q_k} =[a_1,a_2,\dots,a_k]\ , \quad
\frac{r_h}{t_h}=[\xi_1,\xi_2,\dots,\xi_h]\ .$$

Let $\cal Z$ be the generating partition of $X$ in the intervals
$A_0=(0,1/2]$ and $A_1=(1/2,1]$. We consider the AIC of the ergodic
rotation $(X,T,\mu)$, where $\mu$ is the Lebesgue measure, using this
partition.

\begin{theorem} \label{teo:aic_rot}
For all $\alpha \in \R\setminus \Q$ and for $\mu$-almost all $x\in X$
we have $AIC(x,n)\approx$ {\rm constant}.
\end{theorem}
{\em Proof.} By the estimate (\ref{eq:approx_cf}) we have
that the knowledge of the first $k+1$ digits (we remark that
$n$ digits specify $p_i$ and $q_i$ for $i=1,\dots,n$) in the
continued fractions expansion of $\alpha$ allows us to say that
$T^{n}(x)$ is contained in the interval
$$I_n^k:=\left[ (x+n\frac{p_k}{q_k}-\frac{n}{q^2_k a_{k+1}}) ({\rm
mod}\, 1)\, , (x+n\frac{p_k}{q_k}+\frac{n}{q^2_k a_{k+1}}) ({\rm
mod}\, 1) \right]$$ 
The knowledge of the endpoints of this interval yields unambiguously
the symbol $\omega_n$ in the symbolic orbit of $x$ provided $I_n^k$
does not contain $0$ or $1/2$.  Let $T_k$ be the translation of the
angle $\frac{p_k}{q_k}$ and $\omega^k$ be the symbolic orbit of the
point $x$ with the map $T_k$ on $Z$: we have $\omega^k_n
=\omega_n$ provided $I_n^k$ does not contain $0$ or $1/2$.  Set
$$A(k,n) := \left\{ x\in X \, :\, \omega_n\not= \omega^k_n \right\}$$
A short calculation yields
$$\mu (A(k,n)) \le \ \frac{4n}{q^2_k a_{k+1}}$$ 
for all $k$ and $n$. The measure of the set of initial points $x$ that
have at least one symbol among $\omega_1\dots\omega_n$ different from
the corresponding symbols in the word $\omega^k$ is then estimated as
$$\mu \left(\bigcup_{i=1}^n A(k,i)\right) \le \sum_{i=1}^n \mu
(A(k,i)) \le \ \sum_{i=1}^n \frac{4i}{q^2_k
a_{k+1}}=\frac{2n(n+1)}{q^2_k a_{k+1}}$$

Let now $k(n)$ be a sequence of integers such that 
\begin{equation} \label{eq:bc1}
\sum_{n=1}^\infty \frac{4n}{q^2_{k(n)} a_{k(n)+1}} \ <\infty
\end{equation}
By Borel-Cantelli Lemma we have that for $m$-almost all $x\in X$ there
is $n_0>0$ such that $x\in X \setminus A(k,n)$ for all $n\ge n_0$ and
therefore
$$
{AIC}(x,n)={AIC}(\omega_0\omega_1\dots \omega_{n-1}) \le
{AIC}(\omega_o^{k(n)}\omega_1^{k(n)}\dots \omega_{n-1}^{k(n)})+{\rm
const}
$$ 
that is to specify the first $n$ symbols of the sequence $\omega$ it
is enough to specify the first $n$ symbols of the symbolic word
$\omega^{k(n)}$ obtained by the orbit of the point $x$ drawn using the
periodic map $T_{k(n)}$, defined by
$T_{k(n)}(x)=x+\frac{p_{k(n)}}{q_{k(n)}} \, ({\rm mod}\, 1) $.

We have now to estimate the information content of the first $n$
symbols of the word $\omega^k$. First of all we have that to
reconstruct $\omega^k$ it is necessary to know $\frac{p_k}{q_k}$,
hence we need to know the first $k$ symbols of the sequence
$\sigma$. 

Now it remains the uncertainty about the initial condition we use to
generate $\omega^k$. We repeat the same argument as before. Let
$[\xi_1,\xi_2,\dots,\xi_n,\dots]$ be the expansion in continued
fractions of the initial condition $x$. If we know only the first $h$
symbols in this expansion, we can approximate $x$ by the rational
number $\frac{r_h}{t_h}$ so that $x+n\frac{p_k}{q_k} \, ({\rm mod}\,
1)$ is contained in the interval
$$J_n^h := \left[ (\frac{r_h}{t_h} +n\frac{p_k}{q_k} - \frac{1}{t^2_k
\,\xi_{k+1}}) ({\rm mod}\, 1),\ (\frac{r_h}{t_h} + \frac{p_k}{q_k}
+\frac{1}{t^2_k\, \xi_{k+1}}) ({\rm mod}\, 1) \right]$$ Reasoning as
above and choosing $h(n)$ so that
\begin{equation} \label{eq:bc2}
\sum_{n=1}^\infty \frac{4}{t^2_{h(n)} \, \xi_{h(n)+1}} \ <\infty 
\end{equation}
we conclude that for $m$-almost all $x\in X$ there exists $n_1>0$
such that for all $n\ge n_1$ it is possible to determine $\omega^k_n$ 
just knowing $[\xi_1,\xi_2,\dots,\xi_{h(n)}]$. 

From this argument we can conclude that for almost any $x\in X$ we
have that there exist integer functions $k(n)$ and $h(n)$ such that
for $n$ large enough 
\be \label{eq:diseg_aic_rot}
\begin{array}{cl}
{AIC}(x,n) & \le  {AIC}(\omega_0^{k(n)}\omega_1^{k(n)}\dots
\omega_{n-1}^{k(n)})+{\rm const} \le \\[0.4cm]
 & \le AIC(\sigma_{k(n)})+ AIC(s_{h(n)})+ \log_2 n + {\rm const}
\end{array}
\end{equation} 
where $AIC(\sigma_{k(n)})$ and $AIC(s_{h(n)})$ denote the Algorithmic
Information Content of the first $k(n)$ and $h(n)$ symbols of the
sequences $\sigma$ and $s$, respectively, whereas the term $\log_2 n$
is the information needed to reconstruct the $h(n)$-long symbolic
orbit of the rational point $\frac{r_{h(n)}}{t_{h(n)}}$ by the map
$T_{k(n)}$.

We now study the general problem of estimating ${AIC}(\sigma_k)$. We have
\be \label{eq:gen_estim}
{AIC}(\sigma_k) \le k \ \log_2(\max \{ a_i\ : \ i=1,\dots,k \} )
\end{equation}
where we recall that $a_i$ are the partial quotients in the expansion
of $\alpha$. This estimate comes from the general estimate of the
Algorithmic Information Content of a word written on a given
alphabet by the length of the word times the logarithm of the number
of letters in the alphabet.

To conclude the proof we need the following lemma.

\begin{lemma} \label{lem:fraz} Let $f(k)=\max_{i\le k} a_i$ and $k(n)$ a
sequence such that
\begin{equation} \label{eq:cond1}
\sum_{n=1}^\infty \ \frac{4n}{q^2_{k(n)} a_{k(n)+1}} < \infty
\end{equation}
Then it is possible to choose $k(n)$ such that $k(n) \log (f(k(n)))
\approx$ constant.
\end{lemma}

Applying this lemma to $AIC(\sigma_{k(n)})$ and, with an entirely
analogous argument, to $AIC(s_{h(n)})$, we obtain the thesis of the
theorem by equation (\ref{eq:diseg_aic_rot}). \qed

\vsni 
\noindent 
{\em Proof of Lemma \ref{lem:fraz}.} Given $\alpha \in \R
\setminus \Q$, the best bounds we can give on the denominators $q_k$
are \be\label{eq:gauss}
\frac{1}{2\,\prod_{i=0}^{k-1}G^{i}(\alpha)}<q_k<
\frac{1}{\prod_{i=0}^{k-1}G^{i}(\alpha)} 
\end{equation} 
where $G :[0,1]\to [0,1]$ is the Gauss transformation given by $G(x)
=\{1/ x\}$ for $x> 0$ and $G(0)=0$. From this we obtain $q_k \ge
\frac{1}{2} f(k)$ for all $k$. Moreover, as already remarked in
Section 4.1, the golden mean $\bar \alpha = \frac{\sqrt{5}-1}{2}$ is
the number for which the denominators $\bar q_k$ are the smallest and
they are equal to the Fibonacci numbers $F_k$. Hence in general, for
all irrational $\alpha$ it holds $q_k \ge F_k$ for all $k$. 

We have thus obtained two different estimates for the denominators
$q_k$, and the two estimates have the same order if the function
$f(k)$ grows exponentially. 

When $f(k) = {\cal O}(\exp(k))$, we can use the estimate $q_k \ge
F_k$. Hence
$$\sum_{n=1}^\infty \ \frac{4n}{q^2_{k(n)} a_{k(n)+1}} \le
\sum_{n=1}^\infty \ \frac{4n}{F^2_{k(n)}} \sim \sum_{n=1}^\infty \
\frac{n}{\left(\frac{1}{\bar \alpha}\right)^{2 k(n)}}$$ Hence one can
choose $k(n)\sim \log n^{\frac{3}{2}}$ to obtain condition
(\ref{eq:cond1}). Moreover $k \log(f(k)) = {\cal O}(k^2)$, from which
we obtain the thesis of the lemma.

If instead $f(k)$ grows faster than an exponential, then the better
estimate for the denominators $q_k$ is $q_k \ge f(k)$. In this case 
$$\sum_{n=1}^\infty \ \frac{4n}{q^2_{k(n)} a_{k(n)+1}} \le
\sum_{n=1}^\infty \ \frac{4n}{f^2(k(n))}$$ and condition
(\ref{eq:cond1}) is obtained choosing $f(k(n)) \sim
n^{\frac{3}{2}}$. From this is then easy to see that $k(n)\approx$
constant, hence $k(n) \log (f(k(n))) \sim k(n) \log n
\approx$ constant. \qed

\vsni 
From the proof of Theorem \ref{teo:aic_rot} it follows that the
uncertainty on the orbits of points $x\in X$ depends on the
uncertainty on the initial condition $x$, but also on the uncertainty
on the angle $\alpha$, that in principle may need an infinite amount
of information to be known.

If $\alpha$ needs a finite amount of information to be known at any
accuracy we call it {\em constructive} (see Definition
\ref{CONST}). In this case stronger results can be obtained (Section
\ref{gala}).

We now study the more general Interval Exchange Transformations (see
Section \ref{iet}).

To study the behaviour of the AIC for orbits of an IET, we use again
the continued fraction expansions. Given the left end points
$d^1=0,d^2,\dots, d^r$ of the partition $\xi$ of $X$ associated to the
IET, we denote their continued fraction expansions by
$[\delta_{(i,1)},\dots,\delta_{(i,n)},\dots]$, for $i=2,\dots,r$. Let
$\frac{p_{(i,k)}}{q_{(i,k)}}=[\delta_{(i,1)},\dots,\delta_{(i,k)}]$ be
the sequence of rational approximants of $d^i$, then it holds (see
equation (\ref{eq:approx_cf})) for all $i=2,\dots,r$
$$\left| d^i - \frac{p_{(i,k)}}{q_{(i,k)}} \right| \le
\frac{1}{q_{(i,k)}^2 \ \delta_{(i,k+1)}}$$
We denote by $i_0$ the index $i$ such that the approximation
$\frac{1}{q_{(i,k)}^2 \ \delta_{(i,k+1)}}$ is the worst. 

From the definition of an IET, we have that the angles $\alpha_j$ of
rotation for the different intervals can be easily obtained from the
points $d^i$ and the permutation $\sigma$ that defines the
transformation. In particular we find that, using the rational numbers
$\frac{p_{(i,k)}}{q_{(i,k)}}$ instead of the points $d^i$, we have
that the angles $\alpha_j$ are approximated by rational numbers
$\beta_j$, such that for all $j=1,\dots,r$ \be
\label{eq:iet_approx_ang} \left| \alpha_j - \beta_j \right| \le
\frac{r}{q_{(i_0,k)}^2 \ \delta_{(i_0,k+1)}}
\end{equation}

\begin{theorem} \label{teo:iet_aic} For an aperiodic interval exchange
transformation $T$ associated to the permutation $\sigma$ and to the
partition $\xi$, we have that for $\mu$-almost all $x\in X$ (where
$\mu$ denotes the Lebesgue measure) it holds $AIC(x,n)\approx$
{\rm constant}, where the AIC is considered with respect to the partition
$\xi$.
\end{theorem}
{\em Proof.} The proof is analogous to that of Theorem
\ref{teo:aic_rot}. Let $x\in X$ be an initial condition for the IET
and denote by $[\xi_1,\dots,\xi_n,\dots]$ its continued fractions
expansion. Again let $\frac{r_h}{t_h}=[\xi_1,\dots,\xi_h]$ be the
rational approximants of $x$. From equation (\ref{eq:approx_cf}) we
have
$$\left| x - \frac{r_h}{t_h} \right| \le \frac{1}{t_h^2 \ \xi_{h+1}}$$
It is evident that for $\mu$-almost all $x\in X$ there exist $h$ and $k$
such that the rational numbers $\frac{r_h}{t_h}$ and
$\frac{p_{(i,k)}}{q_{(i,k)}}$, for $i=2,\dots,r$, are well ordered
on the real line. Hence we obtain the value of $\omega_0(x)$, that is
the first symbol in the symbolic sequence $\omega(x) \in
\Omega:=\{1,\dots,r\}^{\N}$ associated to the IET by the usual method of
symbolic representation.

Given that we know $\omega_0(x)$, we have that $T(x) \in A^{(k,h)}_1$,
where
$$A^{(k,h)}_1 = \left[ \frac{r_h}{t_h} + \beta_{\omega_0(x)} -
\epsilon_1(h,k),\ \ \frac{r_h}{t_h} + \beta_{\omega_0(x)} +
\epsilon_1(h,k) \right]$$
$$\epsilon_1(h,k)= \frac{1}{t_h^2 \ \xi_{h+1}} +
\frac{r}{q_{(i_0,k)}^2 \ \delta_{(i_0,k+1)}}$$ Hence $\omega_1(x)$ is
well defined if the set $A^{(k,h)}_1$ does not intersect the intervals
$$B(i,k) = \left[ d^i - \frac{1}{q_{(i_0,k)}^2 \ \delta_{(i_0,k+1)}},
\ \ d^i + \frac{1}{q_{(i_0,k)}^2 \ \delta_{(i_0,k+1)}} \right] $$
for all $i=1,\dots, r$. Then if we define
$$P^{(k,h)}_1 = \Big\{ x\in X \ /\ A^{(k,h)}_1 \cap B(i,k) \not=
\emptyset \ \mbox{ for some } i \Big\}$$
we have easily that
\be \label{eq:iet_diseg_1}
\mu (P^{(k,h)}_1) \le r \left( \frac{2}{t_h^2 \ \xi_{h+1}} +
\frac{2(r+1)}{q_{(i_0,k)}^2 \ \delta_{(i_0,k+1)}} \right)
\end{equation}

Let us assume now that we know $\omega_0(x)$ and $\omega_1(x)$. We can
then say that $T^2(x) \in A^{(k,h)}_2$, where now
$$A^{(k,h)}_2 = \left[ \frac{r_h}{t_h} + \beta_{\omega_0(x)} +
\beta_{\omega_1(x)} -
\epsilon_2(h,k),\ \ \frac{r_h}{t_h} + \beta_{\omega_0(x)} +
\beta_{\omega_1(x)} +
\epsilon_2(h,k) \right]$$
$$\epsilon_2(h,k)= \frac{1}{t_h^2 \ \xi_{h+1}} +
\frac{2r}{q_{(i_0,k)}^2 \ \delta_{(i_0,k+1)}}$$ Again $\omega_2(x)$ is
well defined if the set $A^{(k,h)}_2$ does not intersect the intervals
$B(i,k)$. Then we define
$$P^{(k,h)}_2 = \Big\{ x\in X \ /\ A^{(k,h)}_2 \cap B(i,k) \not=
\emptyset \ \mbox{ for some } i \Big\}$$
and obtain
\be \label{eq:iet_diseg_2}
\mu (P^{(k,h)}_2) \le r \left( \frac{2}{t_h^2 \ \xi_{h+1}} +
\frac{2(2r+1)}{q_{(i_0,k)}^2 \ \delta_{(i_0,k+1)}} \right)
\end{equation}

Iterating this argument, we can conclude that for a given $n$, the
measure of the set of initial conditions $x$ for which the knowledge
of $\frac{r_h}{t_h}$ and $\frac{p_{(i,k)}}{q_{(i,k)}}$, for
$i=2,\dots,r$, is not enough to construct
$(\omega_0(x),\dots,\omega_{n-1}(x))$, is estimated from above by \be
\label{eq:iet_diseq_n} \mu \left( \bigcup_{j=1}^n \ P^{(k,h)}_j \right)
\le \sum_{j=1}^n \ r \left( \frac{2}{t_h^2 \ \xi_{h+1}} +
\frac{2(jr+1)}{q_{(i_0,k)}^2 \ \delta_{(i_0,k+1)}} \right) 
\end{equation}

At this point the conclusion of the theorem follows as for Theorem
\ref{teo:aic_rot}. From the Borel-Cantelli Lemma, choosing sequences
$k(n)$ and $h(n)$ such that \be \label{eq:iet_bc} \sum_{n=1}^\infty \
r \left( \frac{2}{t_{h(n)}^2 \ \xi_{h(n)+1}} +
\frac{2(nr+1)}{q_{(i_0,k(n))}^2 \ \delta_{(i_0,k(n)+1)}} \right) <
\infty
\end{equation}
we have that for $\mu$-almost all initial conditions $x$
\be \label{eq:iet_diseq_aic}
\begin{array}{rl}
AIC(x,n,Z)\ \le & \mbox{const }+ \log_2 n + AIC(\xi_1,\dots,\xi_{h(n)}) +
\\[0.3cm]
 & + \sum_{i=2}^r \ AIC(\delta_1^i,\dots,\delta_{k(n)}^i)
\end{array}
\end{equation}
where the constant contains also the information necessary to transmit
the permutation $\sigma$ to which the IET is associated, and the term
$\log_2 n$ accounts for the information necessary to transmit the
number $n$ of iterations that have to be performed using the rational
numbers $\frac{r_h}{t_h}$ and $\frac{p_{(i,k)}}{q_{(i,k)}}$, for
$i=2,\dots,r$.

Simple application of equation (\ref{eq:gen_estim}) and Lemma
\ref{lem:fraz} now conclude the theorem. \qed

\vsni
As said before for the irrational rotations, if all the coefficients
of the IET are constructive numbers, stronger results are
possible. This is shown in the next Section.

\begin{remark} Theorems \ref{teo:aic_rot} and \ref{teo:iet_aic} remain
unchanged if we estimate the AIC with respect to any finite measurable
partition whose sets are finite union of intervals. In the proofs we
just need to add the information coming from the end points of the
intervals of the partition, but Lemma \ref{lem:fraz} applies again for
these points.
\end{remark}

\subsection{Constructive piecewise isometries} \label{gala}

In this section we use Theorem \ref{Ganzo!} and some result from
\cite{Ga} about the relation between information and initial condition sensitivity in the 0-entropy case to estimate the information content of orbits coming from certain constructive piecewise isometries.

Piecewise isometries can be somehow sensitive to initial
conditions because of their discontinuities. The orbit of two nearby
starting points can be separated when the points falls on the
different sides of a discontinuity.

Theorem \ref{Ganzo!} gives a result about quantitative recurrence
rates of measure preserving maps in a neighbourhood of some target
points. The theorem allows to estimate the distance of a typical
orbit from the discontinuity points and then to estimate this kind of
sensitivity.

Now we define the class of maps which will be studied.  A real number
is constructive if it can be approximated at any accuracy by an
algorithm.
\begin{definition}[Constructive numbers] \label{CONST} A number $z\in
{\bf R}$ is said to be constructive if there is an algorithm
$A_z(n):{\bf N}\rightarrow {\bf Q}$ such that $A_z(n)=q $ implies
$|q-z|<2^{-n}$.
\end{definition}

This definition was already given by Turing (\cite{Tur}) and the study
of the computable properties of real numbers was one of the
motivations he had in mind when introducing his famous
computing machines.

\begin{definition} \label{IET}
A constructive piecevise isometry on the interval is a piecewise isometry $T$ with the
following properties: it is Lebesgue measure preserving and both the
discontinuity points $y_1,...,y_n$ and the values $T(y_1),...,T(y_n)$
of $T$ at the discontinuity points are constructive numbers.
\end{definition}

Rotations and interval exchanges whose angles are constructive are
included in the previous definition.

If a map satisfies the above definition, for each constructive point
$x\in [0,1]$ such that $x\notin \{y_1,...,y_n\}$ the value of $T(x)$
can be calculated up to any given accuracy by an algorithm. This
implies the following lemma

\begin{lemma} \label{xprimo} Let $T$ be a constructive piecewise isometry
then there is a constructive number $r\in [0,1]$, such that for each
rational $q$, the orbit of the point $x=q+r \ mod(1) $ $\in [0,1]$ can
be followed at any accuracy by an algorithm, in the following sense:

there is an algorithm (a total recursive function) $A:{\bf Q} \times
{\bf N}\times {\bf Q}\rightarrow \Sigma$ such that $\forall k\in {\bf
N},q\in {\bf Q},\epsilon \in {\bf Q}$ $d(T^k(x),(A(q,k,\epsilon) +r)
mod (1))< \epsilon $.
\end{lemma}
{\em Proof.} For the proof we refer to the results that are in
\cite{Ga} Section 3. Let us consider the following set $B= \{x\in
[0,1],\exists \ i,j \ s.t.\ T^i(x)=y_j\}$ and $X'=[0,1]-B $. Let us
consider a constructive number $r$ such that $r$ is incommensurable
with all $y_i$ and $T(y_i)$. Let us consider the standard
interpretation (\cite{Ga} Def. 3) on $[0,1]$ $I$ given by $ I(s)=\sum 2^{-s_i}$ where
$s=s_1s_2...s_n\in \Sigma $.  Consider the computable structure ${\cal
I} $ containing $I$.  Now let us consider $r$ and $I_r$ given by
\begin{equation} 
I_r(s)=I(s)+r \ (mod\ 1)\label{ierre1}
\end{equation}
then $ I_r$ is a computable interpretation and is contained in $I$.
$I_r(\Sigma) $ is contained in $X'$ and let ${\cal I}_r$ be the
computable structure containing $I_r$.

Now the system $(X',{\cal I}_r, T) $ is constructive (for \cite{Ga}
Definition 9) and $T$ is continuous over $X'$ and the assumptions of
\cite{Ga} Lemma 10 are satisfied.  Then \cite{Ga} Lemma 10 holds for
all $x\in X'$ and then for all $x\in I_r(\Sigma)$, easily implying the statement. \qed

\vsni 
Let us now consider a constructive piecewise isometry as in definition above. Let
$y_1,...,y_l$ be the discontinuity points.  Let us consider two nearby
starting points.  As it was said before, during the iterations of the
map the orbits can be separated by the discontinuity, and the map is
in a certain sense sensitive to initial conditions. This sensitivity
will be estimated by the use of Theorem \ref{Ganzo!}.  Now let us
consider a partition $\alpha $ and suppose that $\alpha $ is finite,
and each set of $\alpha$ is made of a finite union of intervals whose
end points are constructive numbers. Let us consider
${AIC(x,\alpha,n)}$, the information content of the orbit of $x$ with
respect to the partition.

\begin{theorem}
Let $([0,1],T) $ be an IET with flips as in Def. \ref{IET}, then for
Lebesgue a.e. $x\in X$
$$\limsup_{n\to \infty}\frac {AIC(x,\alpha,n)}{\log n}\leq 2.$$
\end{theorem}
{\em Proof.} Refining the partition, the information content
increases, then we can suppose that in the partition $\alpha
=\{[\alpha_1,\alpha _2),[\alpha_2, \alpha _3),...,[\alpha _{q-1},
\alpha_{q}]\} $ each set is smaller than the smallest discontinuity of
$T$.

Let
$$ 
B(n,x,\alpha)=\{ y \in X \ s.t.\ \forall i\leq n, T^i(y)\in
[\alpha_j,\alpha _{j+1}) \Leftrightarrow T^i(x)\in [\alpha_j,\alpha
_{j+1}) \}
$$ 
be the set of points having the same symbolic orbit as $x$ for $n$
steps.

Let us consider two nearest points $x,y$, starting from the same set
of $\alpha $. Since $T$ is a piecewise isometry $x$ and $y$ continue
to have the same symbolic orbit until

1) they are separated by a discontinuity, or

2) they go to different sets of $\alpha$, because a point $\alpha
   _i$ separating different sets of $\alpha$ comes between their
   images ($\exists j,i $ such that $T^j(x)\leq \alpha_i \leq T^j(y)$
   or $T^j(x)\geq \alpha_i \geq T^j(y)$).
  
\noindent Moreover, this implies $ B(n,x,\alpha)$ is an interval for
each $n$.

Since $T$ is a piecewise isometry the two above cases does not happen
for $n$ steps if the orbit of $x$ remains far away enough from the
discontinuity and separation points, i.e. if $\forall i\leq n \
{\min}_j (min (d(T^n(x),y_j),d(T^n(x),\alpha_j)))>d(x,y)$. Now we
apply Theorem \ref{Ganzo!}  to $T$ with the Lebesgue measure (which
has local dimension 1 for each point), where the $x_i$ in the
statement are the union of discontinuity points of $T$ and the
separation points of $\alpha$ ($\{x_i\}=\{y_i\}\cup\{\alpha _i\}$). By
this theorem for a.e. point $x$ the length of the interval $
B(n,x,\alpha)$ decreases eventually slower than $\frac C{n^{p }}$ for
some $C$ and each $p>1$.

Now let us consider $r$ found in Lemma \ref{xprimo}. We remark that
in an interval with length $\ell $ there is at least a rational $r$ of
the kind $r=0.s_1s_2...s_v $ with $v\leq -\log \ell +1$.

By the estimation above about the size of $B(n,x,\alpha)$, for
a.e. $x$, for each $n$ there is $ \tilde x \in B(n,x,\alpha)$ such
that $\tilde x=0.s_1s_2...s_v +r \ ({\rm mod} \ 1)$ is such that $v<p
\log n +C$, $\forall p >1$ and the point $\tilde x$ will have the same
symbolic orbit as $x$.

It is possible to construct a program giving the symbolic orbit of
$x$ with the information contained in the map $T$, in the string
$[s_1s_2...s_v]$ and in the number of steps $n$.

The program giving the symbolic orbit of $x$ will contain a string
$[s_1s_2...s_v]$ ($p \log n +C$ bits, for $p $ near to 1 as we want),
the information that is necessary to know $T$ and $p$ (a constant).

A program which outputs the symbolic sequence associated to the orbit
of $x$ will run as follows:

by Lemma \ref{xprimo} it will follow the orbit of $\tilde x$
with an accuracy $1/k $ obtaining points $\tilde
x_1(\epsilon),...,\tilde x_n(\epsilon)$ following the orbit of $\tilde
x$ at a distance less than $1/k $.  The procedure is repeated
increasing $k$ until $\min_{j,i} (d(\tilde x_i(\epsilon),\alpha_j))>2/k
$ (this assures that the sequence $\tilde x_1(\epsilon),...,\tilde
x_n(\epsilon)$ has the same symbolic orbit as $x$ with respect to the
partition $\alpha $.

Once having the sequence $\tilde x_1(\epsilon),...,\tilde
x_n(\epsilon)$ the symbolic orbit with respect to $\alpha $ is
calculated and this is the output of the computation.
By what is said above the length of the program is then less or equal
to $2\log n+C$. \qed


\begin{thebibliography}{DEGHL}

\bibitem[AB]{AB} {\sc P Alessandri, V Berth\`e}, \, {\it Three
distances theorem and combinatorics on words}, L'Enseignement
Math\'ematique {\bf 44} (1998), 103-132

\bibitem[ACS]{ACS} {\sc V Afraimovich, J R Chazottes, B Saussol}, {\sl
Pointwise dimensions for Poincar\'e recurrences associated with maps
and special flows}, Disc. Cont. Dyn. Syst. - A {\bf 9} (2003), 263-280

\bibitem[BS]{BS} {\sc L Barreira, B Saussol}, \, {\it Hausdorff
dimension of measures via Poincar\'e recurrence}, Commun. Math. Phys.
{\bf 219} (2001), 443-463

\bibitem[BCF]{BCF} {\sc V Berth\`e, N Chekhova, S Ferenczi}, \, {\it
Covering numbers: arithmetics and dynamics for rotations and interval
exchanges}, Journal d'Analyse Math\'ematique {\bf 79} (1999), 1-31

\bibitem[B]{B} {\sc C Bonanno}, {\sl Application of Information
Measures to Chaotic Dynamical Systems}, Ph.D. thesis, University
of Pisa, 2003

\bibitem[BoGa]{BoGa} {\sc C Bonanno, S Galatolo}, {\sl Orbit
complexity for the Manneville maps }, in preparation

\bibitem[Bo]{Bo} {\sc M D Boshernitzan}, {\sl Quantitative recurrence
results}, Invent. Math. {\bf 113} (1993), 617-631

\bibitem[Br]{Br} {\sc A Brudno}, {\sl Entropy and the complexity of
trajectories of dynamical systems}, Russ. Math. Surv. {\bf 2} (1983),
127-151 (English transl.)

\bibitem[CFS]{CFS} {\sc I P Cornfeld, S V Fomin and Ya G Sinai},
\,{\it Ergodic Theory}, Springer Verlag, 1982

\bibitem[Fa]{Fa} {\sc K J Falconer}, {\sl The Geometry of Fractal
Sets}, Princeton University Press, 1981

\bibitem[FHZ]{FHZ} {\sc S Ferenczi, C Holton, L Zamboni}, {\sl
Structure of three-interval exchange transformations III: ergodic and
spectral properties}, preprint (2000)

\bibitem[Fu]{Fu} {\sc H Furstenberg}, {\sl Recurrence in Ergodic
Theory and Combinatorial Number Theory}, Cambridge University Press,
1990

\bibitem[Ga]{Ga} {\sc S Galatolo}, {\sl Complexity, initial condition
sensitivity, dimension and weak chaos in dynamical systems},
Nonlinearity {\bf 16} (2003), 1219-1238

\bibitem[HSV]{HSV} {\sc M Hirata, B Saussol, S Vaienti}, {\sl
Statistics of return times: a general framework and new applications},
Commun. Math. Phys. {\bf 206} (1999), 33-55

\bibitem[KH]{KH} {\sc A Katok, B Hasselblatt}, {\sl Introduction to
the Modern Theory of Dynamical Systems}, Cambridge University Press,
1995

\bibitem[Ka]{Ka} {\sc M Kac}, {\sl On the notion of recurrence in
discrete stochastic processes}, Bull. Amer. Math. Soc. {\bf 53}
(1947), 1002-1010

\bibitem[OW]{OW} {\sc D S Ornstein, B Weiss}, {\sl Entropy and data
compression schemes}, IEEE Trans. Inf. Th. {\bf 39} (1993), 78-83

\bibitem[Pe]{Pe} {\sc K Petersen}, {\sl Ergodic Theory}, Cambridge
University Press, 1983

\bibitem[Pr]{Pr} {\sc T Prellberg}, {\sl Maps of Intervals with
Indifferent Fixed Points: Thermodynamic Formalism and Phase
Transitions}, Ph.D. thesis, Virginia Polytechnic Institute and State
University, 1991

\bibitem[RS]{RS} {\sc A M Rockett, P Sz\"usz}, {\it Continued
Fractions}, World Scientific, 1992

\bibitem[STV]{STV} {\sc B Saussol, S Troubetzkoy, S Vaienti}, {\sl
Recurrence, dimensions and Lyapunov exponents}, J. Stat. Phys. {\bf
106} (2002), 623-634

\bibitem[T]{Tur} {\sc A Turing}, {\sl On computable numbers, with an
application to the ``Entscheidungsproblem''},
Proc. Lond. Math. Soc. {\bf 42} (1936), 230-265

\bibitem[W]{W} {\sc H S White}, {\sl Algorithmic complexity of points
in dynamical systems}, Erg. Th. Dyn. Syst. {\bf 13} (1993), 807-830

\bibitem[Y]{Y} {\sc L-S Young}, {\sl Dimensions, entropy and Lyapunov
exponents}, Erg. Th. Dyn. Syst. {\bf 2} (1982), 109-124

\end{thebibliography}
\end{document}